\documentclass[12pt]{article}
\usepackage{amsmath,amssymb,amsthm,amscd,graphicx}
\usepackage{latexsym}
\usepackage{enumerate}
\usepackage{amsbsy, amsfonts, textcomp}
\usepackage{makeidx}
\newtheorem{theorem}{Theorem}[]
\newtheorem{proposition}[theorem]{Proposition}
\newtheorem{lemma}[theorem]{Lemma}
\newtheorem{corollary}[theorem]{Corollary}

\theoremstyle{definition}

\theoremstyle{remark}
\newtheorem{remark}[theorem]{Remark}

\newcommand{\C}{\mathbb{C}}

\newcommand{\N}{\mathbb{N}}

\def \rank{\operatorname{rank}}

\def \max{\operatorname{max}}

\begin{document}

\begin{center}

\Large

{\bf An inversion algorithm for polynomial maps}

\normalsize
\vspace{1cm}

EL\.ZBIETA ADAMUS \\ Faculty of Applied Mathematics, \\ AGH University of Science and Technology \\
al. Mickiewicza 30, 30-059 Krak\'ow, Poland \\
e-mail: esowa@agh.edu.pl \\

\vspace{0.5cm}
PAWE\L \ BOGDAN \\ Faculty of Mathematics and Computer Science, \\ Jagiellonian University \\
ul. \L ojasiewicza 6, 30-348 Krak\'ow, Poland \\
e-mail: pawel.bogdan@uj.edu.pl \\

\vspace{0.5cm}
TERESA CRESPO \\ Departament d'\`{A}lgebra i Geometria, \\ Universitat de Barcelona \\
Gran Via de les Corts Catalanes 585, 08007 Barcelona, Spain \\
e-mail: teresa.crespo@ub.edu

\vspace{0.5cm}
ZBIGNIEW HAJTO \\ Faculty of Mathematics and Computer Science, \\ Jagiellonian University \\
ul. \L ojasiewicza 6, 30-348 Krak\'ow, Poland \\
e-mail: zbigniew.hajto@uj.edu.pl

\end{center}

\vspace{0.2cm}

\begin{abstract}
We present an algorithmic equivalent statement to the Jacobian conjecture. Given a polynomial map $F$ on an affine space of dimension $n$, our algorithm constructs $n$ sequences of polynomials such that $F$ is invertible if and only if the zero polynomial appears in all $n$ sequences. This algorithm provides then a classification of polynomial automorphisms of affine spaces.

\end{abstract}

\section{Introduction}

The Jacobian Conjecture originated in the question raised by Keller in \cite{K} on the invertibility of polynomial maps with Jacobian determinant equal to~1. The question is still open in spite of the efforts of many mathematicians. We recall in the sequel the precise statement of the Jacobian Conjecture, some reduction theorems and other results we shall use. We refer to \cite{E} for a detailed account of the research on the Jacobian Conjecture and related topics.

Let $K$ be a field and $K[X]=K[X_1,\dots,X_n]$ the polynomial ring in the variables $X_1,\dots,X_n$ over $K$. A \emph{polynomial map} is a map $F=(F_1,\dots,F_n):K^n \rightarrow K^n$ of the form

$$(X_1,\dots,X_n)\mapsto (F_1(X_1,\dots,X_n),\dots,F_n(X_1,\dots,X_n)),$$

\noindent where $F_i \in K[X], 1 \leq i \leq n$. The polynomial map $F$ is \emph{invertible} if there exists a polynomial map $G=(G_1,\dots,G_n):K^n \rightarrow K^n$ such that $X_i=G_i(F_1,\dots,F_n),  1 \leq i \leq n$. We shall call $F$ a \emph{Keller map} if the Jacobian matrix

$$J=\left(\dfrac{\partial F_i}{\partial X_j}\right)_{\substack{1\leq i \leq n \\ 1\leq j \leq n}}$$

\noindent has determinant equal to 1. Clearly an invertible polynomial map $F$ has a Jacobian matrix $J$ with non zero determinant and may be transformed into a Keller map by composition with the linear automorphism with matrix $J(0)^{-1}$.

\vspace{0.5cm}
\noindent {\bf Jacobian Conjecture.} {\it Let $K$ be a field of characteristic zero. A Keller map $F:K^n \rightarrow K^n$ is invertible.}

\vspace{0.5cm}
For $F=(F_1,\dots,F_n) \in K[X]^n$, we define the \emph{degree} of $F$ as $\deg F= \max \{\deg F_i : 1\leq i \leq n\}$. It is known that if $F$ is a polynomial automorphism, then $\deg F^{-1} \leq (\deg F)^{n-1}$ (see \cite{BCW} or \cite{RW}).

The Jacobian conjecture for quadratic maps was proved by Wang in \cite{W}.  We state now the reduction of the Jacobian conjecture to the case of maps of third degree (see \cite{BCW}, \cite{Y}, \cite{D1} and \cite{D2}).

\begin{proposition}\label{red} \begin{enumerate}[a)] \item (Bass-Connell-Wright-Yagzhev) Given a Keller map $F:\C^n \rightarrow \C^n$, there exists a Keller map $\widetilde{F}:\C^N \rightarrow \C^N$, $N\geq n$ of the form $\widetilde{F}=Id+H$, where $H(X)$ is a homogeneous cubic map and having the following property: if $\widetilde{F}$ is invertible, then $F$ is invertible too.
\item (Dru\.zkowski) The cubic part $H$ may be chosen of the form

$$\left( (\sum_{j=1}^N a_{1j} X_j)^3, \dots,(\sum_{j=1}^N a_{Nj} X_j)^3 \right)$$

\noindent and with the matrix $A=(a_{ij})_{\substack{1\leq i \leq N \\ 1\leq j \leq N}}$ satisfying $A^2=0$.

\end{enumerate}

\end{proposition}

Polynomial maps in the Dru\.zkowski form are easier to handle than general cubic homogeneous polynomial maps. However we note the following result.

\begin{proposition}[\cite{E2} Proposition 2.9] Let $r \in \N$. If the Jacobian Conjecture holds for all cubic homogeneous polynomial maps in $r$ variables, then for all $n \in \N$  the Jacobian Conjecture holds for all polynomial maps of the form

$$F=X+(AX)^3$$

\noindent with $A \in \mathrm{M}_n(\C)$ and $\rank A \leq r$.
\end{proposition}

Given a polynomial map $F:\C^n \rightarrow \C^n$, we shall call a polynomial $P \in \C[X_1,\dots,X_n]$ \emph{invariant under $F$} if $P(F_1,\dots,F_n)=P(X_1\dots,X_n)$. A polynomial map $F:\C^n \rightarrow \C^n$ of the form $F=Id+H$ is called a \emph{quasi-translation} if $F^{-1}=Id-H$.

In this paper we present a recursive algorithm to invert polynomial maps. Given a polynomial map $F:\C^n \rightarrow \C^n$ of the form $F=Id+H$, where $H(X)$ is a homogeneous cubic map, our algorithm constructs, for $1\leq i \leq n$, a sequence $P_k^i$ of polynomials in $\C[X]$ with $P_0^i=X_i$ such that $F$ is invertible if and only $P_{m_i}^i=0$, for some integer $m_i>0$, $1\leq i \leq n$ and, when this is the case, we may compute $F^{-1}$ as an alternating sum of the polynomials $P_k^i$, $0\leq k <m_i$. In the last section, we apply the algorithm to several examples of polynomial maps.

\section{The algorithm}

Let us consider a polynomial map $F:\C^n \rightarrow \C^n$.
 Given a polynomial $P(X_1,\dots,X_n) \in \C[X]=\C[X_1,\dots,X_n]$, we define the following sequence of polynomials in $\C[X]$,

$$\begin{array}{lll} P_0(X_1,\dots,X_n) &= &P(X_1,\dots,X_n), \\
P_1(X_1,\dots,X_n) &=& P_0(F_1,\dots,F_n)-P_0(X_1,\dots,X_n), \end{array}$$

\noindent and, assuming $P_{k-1}$ is defined,

$$P_k(X_1,\dots,X_n) = P_{k-1}(F_1,\dots,F_n)-P_{k-1}(X_1,\dots,X_n).$$

The following lemma is easy to prove.

\begin{lemma}\label{lem} For a positive integer $m$, we have

$$P(X_1,\dots,X_n)= \sum_{l=0}^{m-1} (-1)^l P_l(F_1,\dots,F_n)+(-1)^m P_m(X_1,\dots,X_n).$$

\noindent In particular, if we assume that for some integer $m$, $P_m(X_1,\dots,X_n)=0$, then

$$P(X_1,\dots,X_n)= \sum_{l=0}^{m-1} (-1)^l P_l(F_1,\dots,F_n).$$

\end{lemma}

\begin{theorem}\label{theo} Let $F:\C^n \rightarrow \C^n$ be a polynomial map of the form

$$\left\{ \begin{array}{lll} F_1(X_1,\dots,X_n)&=& X_1+H_1(X_1,\dots,X_n) \\ F_2(X_1,\dots,X_n)&=& X_2+H_2(X_1,\dots,X_n)\\ & \vdots & \\  F_n(X_1,\dots,X_n)&=& X_n+H_n(X_1,\dots,X_n) , \end{array} \right.$$

\noindent where $H_i(X_1,\dots,X_n)$ is a homogeneous polynomial in $X_1,\dots,X_n$ of degree~$3$, $1\leq i \leq n$. Let us assume that $F$ is invertible. Then for the polynomial sequence $(P_k^i)$ constructed with $P=X_i$, there exists an integer $m_i$ such that $P_{m_i}^i = 0$ and the inverse map $G$ of $F$ is given by

$$G_i(Y_1,Y_2,\dots,Y_n)= \sum_{l=0}^{m_i-1} (-1)^l P_l^i(Y_1,Y_2,\dots,Y_n), \, 1 \leq i \leq n.$$

\end{theorem}

\noindent {\it Proof.} Taking into account lemma \ref{lem}, it is enough to prove that there exists an integer $m_i$ such that $P_{m_i}^i=0$, for $i=1,\dots,n$. If $F$ is invertible, then $G=F^{-1}$ is of the form

$$\left\{ \begin{array}{lll} G_1(Y_1,Y_2,\dots,Y_n)&=& Y_1+J_1 \\ G_2(Y_1,Y_2,\dots,Y_n)&=& Y_2+J_2 \\ & \vdots & \\  G_n(Y_1,Y_2,\dots,Y_n) &=& Y_n+J_n , \end{array} \right.$$

\noindent where $J_i$ is a polynomial in $Y_1,\dots,Y_n$, $1\leq i \leq n$. Let $N_i=\deg J_i$. Let us consider, for a fixed $i$, the polynomial sequence

$$\begin{array}{lll} P_0^i(X_1,\dots,X_n) &= & X_i, \\
P_1^i(X_1,\dots,X_n) &=& F_i(X_1,\dots,X_n)-X_i=H_i(X_1,\dots,X_n), \\
P_2^i(X_1,\dots,X_n) &=& H_i(F_1,\dots,F_n)-H_i(X_1,\dots,X_n), \\
&\vdots &
\end{array}$$

\noindent We write the Taylor series for the polynomial $H_i(F_1,\dots,F_n)=H_i(X_1+H_1,\dots,X_n+H_n)$ and obtain

$$\begin{array}{lll}P_2(X_1,\dots,X_n) &=& H_i(F_1,\dots,F_n)-H_i(X_1,\dots,X_n) \\ &=& Q_{21}+Q_{22}+Q_{23}
\end{array}$$

\noindent where

$$\begin{array}{lll} Q_{21} &=& \sum_{j=1}^n H_j \dfrac{\partial H_i}{\partial X_j} \\ [10pt] Q_{22} &=& \dfrac 1 2 \sum_{1\leq j,k \leq n} H_j H_k \dfrac{\partial^2 H_i}{\partial X_j\partial X_k} \\ [10pt] Q_{23} &=& \dfrac 1 6 \sum_{1\leq j,k,l \leq n} H_jH_kH_l \dfrac{\partial^3 H_i}{\partial X_j \partial X_k \partial X_l} \end{array} $$

\noindent and the polynomials $Q_{21},Q_{22},Q_{23}$ are either zero or homogeneous polynomials of degrees $5,7,9$, respectively. Let us prove by induction that

$$P_n^i= \sum_{j=1}^{(3^n-2n+1)/2} Q_{nj},$$

\noindent where $Q_{nj}$ is either zero or a homogeneous polynomial of degree $2n-1+2j$. We have already seen it for $n=2$. Let us assume $P_{n-1}^i= \sum_{j=1}^{(3^{n-1}-2n+3)/2} Q_{n-1,j}$ with $Q_{n-1,j}$  either zero or a homogeneous polynomial of degree $2n-3+2j$. We want to prove the property for $P_n^i$. We have

$$\begin{array}{l} P_n^i(X_1,\dots,X_n) = P_{n-1}^i(F_1,\dots,F_n)-P_{n-1}^i(X_1,\dots,X_n) \\ [8pt] = \sum_{j=1}^{(3^{n-1}-2n-1)/2} Q_{n-1,j}(F_1,\dots,F_n)- Q_{n-1,j}(X_1,\dots,X_n). \end{array}$$

\noindent If $Q(X_1,\dots,X_n)$ is a homogeneous polynomial of degree $d$,

$$Q(F_1,\dots,F_n)-Q(X_1,\dots,X_n)= \sum_{j=1}^n H_j \dfrac{\partial Q}{\partial X_j}+ \dots$$

\noindent is a sum of polynomials which are either zero or homogeneous polynomials of degrees $3+d-1, 6+d-2, \dots, 3d$. Hence $Q_{n-1,j}(F_1,\dots,F_n)- Q_{n-1,j}(X_1,\dots,X_n)$ is a sum of polynomials which are either zero or homogeneous polynomials of degrees $2n+2j-1, 2n+2j+1, \dots, 6n-9+6j$. Therefore $P_n^i= \sum_{j=1}^{(3^n-2n+1)/2} Q_{nj}$ where $Q_{nj}$ is either zero or a homogeneous polynomial of degree $2n-1+2j.$

Applying lemma \ref{lem}, we obtain

$$X_i= \sum_{l=0}^{m-1} (-1)^l P_l(F_1,\dots,F_n)+(-1)^m P_m(X_1,\dots,X_n).$$

\noindent Now, we have that $X_i=G_i(F_1,\dots,F_n)$ and $G_i$ is a polynomial of degree $N_i$. Then the summands appearing in the development of the expression \linebreak $G_i(F_1(X_1,\dots,X_n),\dots,F_n(X_1,\dots,X_n))$ have total degree $3N_i$ in $X_1,\dots,X_n$. If $m > (3N_i-1)/2$, then $P_m$ is a sum of polynomials which are either zero or homogeneous of degrees bigger than $3N_i$. Hence $P_m=0$.
$\Box$

\vspace{0.5cm}
We shall give now a bound for the degrees of the polynomials $P_k^i$ in Theorem \ref{theo}.

\begin{proposition} Let $F, G$ and $(P_k^i)$ be as in Theorem \ref{theo}. Let $N_i=\deg G_i$. Then $\deg P_k^i \leq 3N_i$. In particular, if $P_{(3N_i-1)/2} \neq 0$, then it is homogeneous of degree $3N_i$.
\end{proposition}

\noindent {\it Proof.} We have seen in the proof of Theorem \ref{theo} that $\deg P_k^i \leq 3N_i$, for $k\leq [log_3 N_i]+1$. Let $m=[log_3 N_i]+2$. Applying Lemma \ref{lem}, we obtain

$$X_i= \sum_{l=0}^{m-1} (-1)^l P_l^i(F_1,F_2,\dots,F_n)+(-1)^m P_m^i(X_1,\dots,X_n).$$

\noindent Reasoning as in the proof of Theorem \ref{theo}, we obtain that the homogeneous summands of $P_m^i$ of degrees $>3 N_i$ must be zero, hence $\deg P_m^i \leq 3N_i$. By induction we obtain $\deg P_k^i \leq 3N_i$ for all $k$. By Theorem \ref{theo}, the homogeneous summands of $P_{(3N_i-1)/2}$ have degrees $\geq 3N_i$, so we obtain the last sentence of the Proposition. $\Box$

\begin{corollary} Let $F$ be a polynomial map as in Theorem \ref{theo}. If $F$ is invertible, then there exists $P \in \C[X_1,\dots,X_n]$ invariant under $F$.
\end{corollary}

\noindent {\it Proof.} Indeed, if $P_m^i$ is the first zero polynomial in the sequence $(P_k^i)$, then $P_{m-1}^i$ is invariant under $F$. $\Box$

\begin{remark} We have presented our algorithm for cubic homogeneous complex polynomial maps but one can see easily that it is valid for cubic homogeneous polynomial maps defined over any field $K$ of characteristic 0.
\end{remark}

Due to Proposition \ref{red}, we have presented the algorithm in the cubic homogeneous case however its validity is wider.

\begin{proposition}\label{gen} Let $F:\C^n \rightarrow \C^n$ be a polynomial map of the form

$$\left\{ \begin{array}{lll} F_1(X_1,\dots,X_n)&=& X_1+H_1(X_1,\dots,X_n) \\ F_2(X_1,\dots,X_n)&=& X_2+H_2(X_1,\dots,X_n)\\ & \vdots & \\  F_n(X_1,\dots,X_n)&=& X_n+H_n(X_1,\dots,X_n) , \end{array} \right.$$

\noindent where $H_i(X_1,\dots,X_n)$ is a polynomial in $X_1,\dots,X_n$ of lower degree $d_i \geq 2$, $1\leq i \leq n$. Let us assume that $F$ is invertible. Then for the polynomial sequence $(P_k^i)$ constructed with $P=X_i$, there exists an integer $m_i$ such that $P_{m_i}^i = 0$ and the inverse map $G$ of $F$ is given by

$$G_i(Y_1,Y_2,\dots,Y_n)= \sum_{l=0}^{m_i-1} (-1)^l P_l^i(Y_1,Y_2,\dots,Y_n), \, 1 \leq i \leq n.$$

\end{proposition}

\noindent {\it Proof.} Following the proof of Theorem \ref{theo}, one can see that the lower degrees of the polynomials  $(P_k^i)$ are increasing by (at least) $d-1 >0$, where $d=\min \{d_i\}$, hence we obtain $P_m^i=0$, for $m>(DN_i-d_i)/(d-1)+1$, where $D=\deg F, N_i=\deg G_i$. $\Box$

\vspace{0.5cm}
As a consequence of Proposition \ref{gen}, we obtain a direct proof of the characterization of quasi-translations obtained by de Bondt in \cite{dB}, Proposition 1.1. See also \cite{dB2} Proposition 1.3.

\begin{proposition}\label{quasi} Let $F:\C^n \rightarrow \C^n$ be a polynomial map of the form $F=Id +H$, where $H_i(X_1,\dots,X_n)$ is a polynomial in $X_1,\dots,X_n$ of lower degree $d_i \geq 2$, $1\leq i \leq n$. Then the following conditions are equivalent.

\begin{enumerate}[a)]
\item $F$ is a quasi-translation.
\item $JH\cdot H=0$, where $JH$ denotes the jacobian matrix of $H$.
\item $P_2^i=0$, $1\leq i \leq n$, where $(P_k^i)$ denotes the polynomial sequence constructed with $P=X_i$.
\end{enumerate}
\end{proposition}

\noindent {\it Proof.} From Lemma \ref{lem} applied to $P(X)=X_i$, we have $X_i=F_i-H_i(F)+P_2^i(X)$, hence the equivalence between a) and c) is clear. Now

$$P_2^i(X)=H_i(F)-H_i(X)=\sum_{k\geq 1}\sum_{1\leq j_1,\dots,j_k \leq n} \dfrac {\partial^k H_i}{\partial X_{j_1} \dots \partial X_{j_k}} H_{j_1} \dots H_{j_k}.$$

\noindent Let us assume b). We have then $\sum_{j=1}^n \dfrac {\partial H_i}{\partial X_{j}}H_j=0$, $1\leq i \leq n$. We shall prove by induction

\begin{equation}\label{eq} \sum_{1\leq j_1,\dots,j_k \leq n} \dfrac {\partial^k H_i}{\partial X_{j_1} \dots \partial X_{j_k}}H_{j_1} \dots H_{j_k}=0
\end{equation}

\noindent for all $k\geq 0$. For $k=1$, it is true by hypothesis. Let us assume that it is true for $k$. By applying $\partial /\partial X_j$ to (\ref{eq}), we obtain

$$\sum_{1\leq j_1,\dots,j_k \leq n} (\dfrac {\partial^{k+1} H_i}{\partial X_{j_1} \dots \partial X_{j_k}\partial X_j}H_{j_1} \dots H_{j_k}+\sum_{l=1}^k \dfrac {\partial^k H_i}{\partial X_{j_1} \dots \partial X_{j_k}}H_{j_1}\dots \dfrac {\partial H_{j_l}}{\partial X_j} \dots H_{j_k})=0.$$

\noindent Multiplying by $H_j$, summing from $j=1$ to $n$ and taking into account the case $k=1$ with $i=j_1,\dots,j_k$, we obtain (\ref{eq}) for $k+1$.

Let us assume c). If $H$ is homogeneous, the summands of $P_2^i$ are homogeneous of different degrees, hence all of them vanish and since the first summand $\sum_{j=1}^n \dfrac {\partial H_i}{\partial X_{j}}H_j$ is the $i$-th row of $JH\cdot H$, we obtain b). If $H$ is not homogeneous, we obtain successively the vanishing of the homogeneous summands of $\sum_{j=1}^n \dfrac {\partial H_i}{\partial X_{j}}H_j$ in increasing degree order. $\Box$

\begin{remark} In the set of polynomial maps as in Proposition \ref{quasi}, our algorithm provides a filtration $\{ \mathcal{P}_n \}_{n\geq 1}$ such that $\mathcal{P}_n$ is the subset of polynomial maps with $P_n^i=0, 1\leq i \leq n$. Clearly $\mathcal{P}_1 = \{ Id \}$ and, applying Proposition \ref{quasi}, $\mathcal{P}_2$ is the set of quasi-translations.
\end{remark}

\section{Examples}

\subsection{} We consider the following Keller map in dimension 2 to illustrate that the validity of our algorithm is not restricted to homogeneous cubic maps.

$$\left\{ \begin{array}{lll} F_1 &= & X_1+(X_2+X_1^3)^2 \\ F_2 &=& X_2 + X_1^3 \end{array} \right.$$

\noindent Let us write $H_1:= (X_2+X_1^3)^2, H_2:=X_1^3$. Taking into account the bound for the degree of the inverse, we may disregard the terms of degree $> 6$ in the polynomial sequences $(P_k^i)$. We obtain

$$\begin{array}{lll} P_0^1&=& X_1 \\ P_1^1 &=& H_1 = (X_2+X_1^3)^2 \\ P_2^1 &=& 3X_1^6+6X_1X_2^5+6X_1^2X_2^3+2X_1^3X_2 \\
 P_3^1 &=& 2X_1^6+18X_1X_2^5+6X_1^2X_2^3 \\ P_4^1 &=& 12 X_1 X_2^5 \\ P_5^1 &=& 0 \end{array}$$

\noindent which gives $X_1= P_0^1(F_1,F_2)-P_1^1(F_1,F_2)+P_2^1(F_1,F_2)-P_3^1(F_1,F_2)+P_4^1(F_1,F_2)= Y_1-Y_2^2$ and

$$\begin{array}{lll} P_0^2&=& X_2 \\ P_1^2 &=& H_2 = X_1^3\\ P_2^2 &=& 6X_1^5X_2+X_2^6+3X_1X_2^4+3X_1^2X_2^2 \\ P_3^2 &=& 6X_1^5 X_2+6 X_2^6+6X_1X_2^4 \\ P_4^2 &=& 6X_2^6 \\ P_5^2 &=& 0 \end{array}$$

\noindent which gives $X_2= P_0^2(F_1,F_2)-P_1^2(F_1,F_2)+P_2^2(F_1,F_2)-P_3^2(F_1,F_2)+P_4^2(F_1,F_2)= Y_2-Y_1^3+Y_2^6-3Y_1Y_2^4+3Y_1^2Y_2^2=Y_2-(Y_1-Y_2^2)^3.$

\subsection{} We consider the following Keller map $F$ in dimension 5.

$$\left\{ \begin{array}{lll} F_1 &= & X_1+a_1X_4^3+a_2x_4^2X_5+a_3X_4X_5^2+a_4x_5^3 +\dfrac{a_2c_5X_2X_4^2}{c_2}+\dfrac {2a_3c_5X_2X_4X_5}{c_2} \\[10pt] && +\dfrac{3a_4c_5X_2X_5^2}{c_2}-\dfrac{a_2e_2X_3X_4^2}{c_2}-\dfrac{2a_3e_2X_3X_4X_5}{c_2}-\dfrac{3a_4e_2X_3X_5^2}{c_2}\\[10pt] && +\dfrac {a_3c_5^2X_2^2X_4} {c_2^2}+\dfrac {3a_4c_5^2X_2^2X_5} {c_2^2}-\dfrac {2a_3c_5e_2X_2X_3X_4} {c_2^2}-\dfrac {6a_4c_5e_2X_2X_3X_5} {c_2^2} \\[10pt] && +\dfrac {a_3e_2^2X_3^2X_4} {c_2^2}+\dfrac {3a_4e_2^2X_3^3X_5} {c_2^2}+\dfrac {a_4e_5^3X_2^3}{c_2^3}-\dfrac {3a_4c_5^2e_2X_2^2X_3}{c_2^3}\\[10pt] &&+\dfrac {3a_4c_5e_2^2X_2X_3^2}{c_2^3}-\dfrac {a_4e_2^3X_3^3}{c_2^3} \\ [10pt] F_2 &=& X_2 + b_1X_4^3 \\ [10pt] F_3 &=& X_3+c_5X_2X_4^2+c_1X_4^3+c_2X_4^4X_5-e_2X_3X_4^2 \\ [10pt] F_4 &=& X_4 \\ [10pt] F_5 &=& X_5+e_2X_4^2X_5 +\dfrac {c_5e_2X_2X_4^2}{c_2} -\dfrac {e_2^2X_3X_4^2}{c_2}- \dfrac {(b_1c_5-c_1e_2)X_4^3}{c_2} \end{array} \right.$$

\noindent with parameters $a_1,a_2,a_3,a_4, b_1, c_1, c_2, c_5, e_2$. By applying the algorithm we obtain $P_2^i=0$, for $i=1,\dots,5$, hence $F$ is a quasi-translation.

\subsection{} We consider the following Keller map $F$ in dimension 6

$$\left\{ \begin{array}{lll}
F_1 &=& X_1+a_5e_1(X_1+X_2)^3/a_4+a_4X_2X_4X_6+a_5X_4X_5X_6 \\
F_2 &=& X_2-a_5e_1(X_1+X_2)^3/a_ 4 \\
F_3 &=& X_3+c_1X_1^3+c_2(X_1+X_5)^3+c_3(X_1+X_2)^3+c_4(X_1+X_4)^3+c_5X_6^3 \\
F_4 &=& X_4+d_4X_2X_6^2+a_5d_4X_5X_6^2/a_4 \\
F_5 &=& X_5+e_1(X_1+X_2)^3 \\
F_6 &=& X_6
\end{array} \right.$$

\noindent with parameters $a_4,a_5,c_1,c_2,c_3,c_4,c_5,d_4,e_1$. Denoting $G=F^{-1}$ and taking variables $(Y_1,\dots,Y_6)$ for $G$, we obtain $P_8^1=0$, $P_9^2=0$ and

\footnotesize

$$\begin{array}{l}G_1 = -(20a_5^4e_1 Y_6^9d_4^3Y_2^3Y_5^3a_4^3-6a_5^4e_1Y_6^4Y_4Y_5^3d_4Y_2a_4^2+a_5e_1a_4^6Y_6^9Y_2^6d_4^3+3Y_6^3a_5e_1Y_2^4d_4a_4^4\\  +3a_5e_1a_4^5Y_6^6Y_2^5d_4^2-a_5e_1a_4^6Y_2^3Y_4^3Y_6^3+a_4^5Y_2Y_4Y_6-Y_6^3a_4^5d_4Y_2^2+a_5^7e_1Y_6^9d_4^3Y_5^6 \\  +a_5e_1a_4^3Y_1^3+a_5e_1a_4^3Y_2^3-Y_1a_4^4-3a_5^2e_1a_4^5Y_2^2Y_4^3Y_6^3Y_5-6a_5e_1Y_1Y_2^2Y_4Y_6a_4^4 \\  +3a_5^5e_1Y_6^6Y_1d_4^2Y_5^4a_4+12a_5^4e_1Y_6^6Y_2^2d_4^2Y_5^3a_4^2-6a_5^4e_1Y_6^4Y_1Y_4d_4Y_5^3a_4^2+3a_5e_1a_4^3Y_1Y_2^2 \\ -
6a_5e_1a_4^5Y_6^4Y_1Y_2^3Y_4d_4+3Y_6^3a_5e_1Y_1^2Y_2^2d_4a_4^4+3a_5e_1a_4^5Y_1Y_2^2Y_4^2Y_6^2-30a_5^3e_1Y_6^7Y_4d_4^2Y_2^3Y_5^2a_4^4 \\  +6Y_6^3a_5e_1Y_1Y_2^3d_4a_4^4+
15a_5^3e_1Y_6^9d_4^3Y_2^4Y_5^2a_4^4+6a_5^2e_1a_4^5Y_6^9Y_2^5d_4^3Y_5-18a_5^2e_1Y_6^4Y_2^3Y_4d_4Y_5a_4^4\\  +12a_5^2e_1a_4^5Y_6^5Y_2^3Y_4^2d_4Y_5-
15a_5^2e_1a_4^5Y_6^7Y_2^4Y_4d_4^2Y_5+12a_5^2e_1Y_6^6Y_1Y_2^3d_4^2Y_5a_4^4-18a_5^2e_1Y_6^4Y_1Y_2^2Y_4d_4Y_5a_4^4 \\  -3a_5e_1a_4^6Y_6^7Y_2^5Y_4d_4^2+
3a_5e_1a_4^3Y_1^2Y_2+18a_5^3e_1Y_6^5Y_4^2d_4Y_2^2Y_5^2a_4^4+12a_5^2e_1Y_6^6Y_2^4d_4^2Y_5a_4^4\\  +3a_5e_1a_4^5Y_6^6Y_1Y_2^4d_4^2+3a_5e_1a_4^6Y_6^5Y_2^4Y_4^2d_4-
6a_5e_1Y_6^4a_4^5Y_2^4Y_4d_4+a_5Y_4Y_5Y_6a_4^4\\  +15a_5^5e_1Y_6^9d_4^3Y_2^2Y_5^4a_4^2+6a_5^6e_1Y_6^9Y_2d_4^3Y_5^5a_4-18a_5^3e_1Y_6^4Y_2^2Y_4d_4Y_5^2a_4^3+
18a_5^3e_1Y_6^6Y_1d_4^2Y_2^2Y_5^2a_4^3\\  -3a_5^6e_1Y_6^7Y_4Y_5^5d_4^2a_4-18a_5^3e_1Y_6^4Y_1Y_4d_4Y_2Y_5^2a_4^3+6a_5^3e_1Y_6^3Y_1d_4Y_5^2Y_2a_4^2-
30a_5^4e_1Y_6^7Y_4d_4^2Y_2^2Y_5^3a_4^3 \\  +12a_5^4e_1Y_6^5Y_2Y_4^2d_4Y_5^3a_4^3+12a_5^4e_1Y_6^6Y_1d_4^2Y_2Y_5^3a_4^2+3a_5^5e_1Y_6^6d_4^2Y_5^4Y_2a_4+
3a_5^3e_1Y_6^3Y_1^2d_4Y_5^2a_4^2 \\  +18a_5^3e_1Y_6^6Y_2^3d_4^2Y_5^2a_4^3+6a_5^2e_1Y_1Y_2Y_4^2Y_6^2Y_5a_4^4-3a_5e_1Y_2^3Y_4Y_6a_4^4-2Y_6^3Y_2a_5d_4Y_5a_4^4 \\  +
3a_5^5e_1Y_6^5Y_4^2Y_5^4d_4a_4^2-3a_5^3e_1Y_2Y_4^3Y_6^3Y_5^2a_4^4-3a_5e_1Y_1^2Y_2Y_4Y_6a_4^4+3a_5e_1a_4^5Y_2^3Y_4^2Y_6^2 \\  +6a_5^2e_1Y_2^2Y_4^2Y_6^2Y_5a_4^4-
3a_5^2e_1Y_4Y_5Y_6Y_2^2a_4^3-a_5^2Y_6^3d_4Y_5^2a_4^3-3a_5^2e_1Y_1^2Y_4Y_5Y_6a_4^3 \\  +3a_5^3e_1Y_4^2Y_5^2Y_6^2Y_2a_4^3+6Y_6^3a_5^2e_1Y_1^2Y_2d_4Y_5a_4^3-
15a_5^5e_1Y_6^7Y_2Y_4d_4^2Y_5^4a_4^2+12Y_6^3a_5^2e_1Y_1Y_2^2d_4Y_5a_4^3 \\  +6Y_6^3a_5^2e_1Y_2^3d_4Y_5a_4^3+3a_5^3e_1Y_6^3d_4Y_5^2Y_2^2a_4^2+
3a_5^3e_1Y_1Y_4^2Y_5^2Y_6^2a_4^3-6a_5^2e_1Y_1Y_4Y_5Y_6Y_2a_4^3 \\  -a_5^4e_1Y_4^3Y_5^3Y_6^3a_4^3)/a_4^4; \end{array}$$

$$\begin{array}{l}G_2=(20a_5^4e_1Y_6^9d_4^3Y_2^3Y_5^3a_4^3-6a_5^4e_1Y_6^4Y_4Y_5^3d_4Y_2a_4^2+a_5e_1a_4^6Y_6^9Y_2^6d_4^3+3Y_6^3a_5e_1Y_2^4d_4a_4^4\\
+3a_5e_1a_4^5Y_6^6Y_2^5d_4^2-
a_5e_1a_4^6Y_2^3Y_4^3Y_6^3+a_5^7e_1Y_6^9d_4^3Y_5^6+a_5e_1a_4^3Y_1^3\\ +a_5e_1a_4^3Y_2^3-3a_5^2e_1a_4^5Y_2^2Y_4^3Y_6^3Y_5-6a_5e_1Y_1Y_2^2Y_4Y_6a_4^4+
3a_5^5e_1Y_6^6Y_1d_4^2Y_5^4a_4\\ +12a_5^4e_1Y_6^6Y_2^2d_4^2Y_5^3a_4^2-6a_5^4e_1Y_6^4Y_1Y_4d_4Y_5^3a_4^2+3a_5e_1a_4^3Y_1Y_2^2-6a_5e_1a_4^5Y_6^4Y_1Y_2^3Y_4d_4 \\ +
3Y_6^3a_5e_1Y_1^2Y_2^2d_4a_4^4+3a_5e_1a_4^5Y_1Y_2^2Y_4^2Y_6^2-30a_5^3e_1Y_6^7Y_4d_4^2Y_2^3Y_5^2a_4^4+6Y_6^3a_5e_1Y_1Y_2^3d_4a_4^4 \\+
15a_5^3e_1Y_6^9d_4^3Y_2^4Y_5^2a_4^4+6a_5^2e_1a_4^5Y_6^9Y_2^5d_4^3Y_5-18a_5^2e_1Y_6^4Y_2^3Y_4d_4Y_5a_4^4+12a_5^2e_1a_4^5Y_6^5Y_2^3Y_4^2d_4Y_5 \\-
15a_5^2e_1a_4^5Y_6^7Y_2^4Y_4d_4^2Y_5+12a_5^2e_1Y_6^6Y_1Y_2^3d_4^2Y_5a_4^4-18a_5^2e_1Y_6^4Y_1Y_2^2Y_4d_4Y_5a_4^4-3a_5e_1a_4^6Y_6^7Y_2^5Y_4d_4^2 \\+
3a_5e_1a_4^3Y_1^2Y_2+18a_5^3e_1Y_6^5Y_4^2d_4Y_2^2Y_5^2a_4^4+12a_5^2e_1Y_6^6Y_2^4d_4^2Y_5a_4^4+3a_5e_1a_4^5Y_6^6Y_1Y_2^4d_4^2\\ +3a_5e_1a_4^6Y_6^5Y_2^4Y_4^2d_4-
6a_5e_1Y_6^4a_4^5Y_2^4Y_4d_4+15a_5^5e_1Y_6^9d_4^3Y_2^2Y_5^4a_4^2+6a_5^6e_1Y_6^9Y_2d_4^3Y_5^5a_4\\ -18a_5^3e_1Y_6^4Y_2^2Y_4d_4Y_5^2a_4^3+
18a_5^3e_1Y_6^6Y_1d_4^2Y_2^2Y_5^2a_4^3-3a_5^6e_1Y_6^7Y_4Y_5^5d_4^2a_4-18a_5^3e_1Y_6^4Y_1Y_4d_4Y_2Y_5^2a_4^3\\ +6a_5^3e_1Y_6^3Y_1d_4Y_5^2Y_2a_4^2-
30a_5^4e_1Y_6^7Y_4d_4^2Y_2^2Y_5^3a_4^3+12a_5^4e_1Y_6^5Y_2Y_4^2d_4Y_5^3a_4^3+12a_5^4e_1Y_6^6Y_1d_4^2Y_2Y_5^3a_4^2\\ +3a_5^5e_1Y_6^6d_4^2Y_5^4Y_2a_4+
3a_5^3e_1Y_6^3Y_1^2d_4Y_5^2a_4^2+18a_5^3e_1Y_6^6Y_2^3d_4^2Y_5^2a_4^3+6a_5^2e_1Y_1Y_2Y_4^2Y_6^2Y_5a_4^4\\ -3a_5e_1Y_2^3Y_4Y_6a_4^4+3a_5^5e_1Y_6^5Y_4^2Y_5^4d_4a_4^2-
3a_5^3e_1Y_2Y_4^3Y_6^3Y_5^2a_4^4-3a_5e_1Y_1^2Y_2Y_4Y_6a_4^4\\ +3a_5e_1a_4^5Y_2^3Y_4^2Y_6^2+6a_5^2e_1Y_2^2Y_4^2Y_6^2Y_5a_4^4-3a_5^2e_1Y_4Y_5Y_6Y_2^2a_4^3-
3a_5^2e_1Y_1^2Y_4Y_5Y_6a_4^3\\ +3a_5^3e_1Y_4^2Y_5^2Y_6^2Y_2a_4^3+6Y_6^3a_5^2e_1Y_1^2Y_2d_4Y_5a_4^3-15a_5^5e_1Y_6^7Y_2Y_4d_4^2Y_5^4a_4^2+
12Y_6^3a_5^2e_1Y_1Y_2^2d_4Y_5a_4^3\\ +6Y_6^3a_5^2e_1Y_2^3d_4Y_5a_4^3+3a_5^3e_1Y_6^3d_4Y_5^2Y_2^2a_4^2+3a_5^3e_1Y_1Y_4^2Y_5^2Y_6^2a_4^3-6a_5^2e_1Y_1Y_4Y_5Y_6Y_2a_4^3\\ -
a_5^4e_1Y_4^3Y_5^3Y_6^3a_4^3+Y_2a_4^4)/a_4^4.\end{array}$$

\normalsize

Now,

$$\begin{array}{lll} G_6 &=& Y_6 \\ G_5 &=& Y_5-e_1(G_1+G_2)^3 \\ G_4 &=& Y_4-d_4G_2G_6^2-a_5d_4G_5G_6^2/a_4 \\ G_3 &=& Y_3-c_1G_1^3-c_2 (G_1+G_5)^3-c_3(G_1+G_2)^3-c_4(G_1+G_4)^3+c_5G_6^3 .\end{array}$$

\vspace{0.5cm}
\noindent
{\bf Acknowledgments.} E. Adamus acknowledges support of the Polish Ministry of Science and Higher Education. T. Crespo and Z. Hajto acknowledge support of grant MTM2012-33830, Spanish Science Ministry.

\end{document}